\documentclass[11pt]{amsart}
\renewcommand{\bar}{\overline}
\usepackage{amsmath,amsthm,amscd,euscript}
\def \r{\mathbb R}
\def \q{\mathbb Q}
\def \z{\mathbb Z}
\def \n{\mathbb N}

\newtheorem{theorem}{Theorem}[section]
\newtheorem{lemma}[theorem]{Lemma}

\newtheorem{proposition}[theorem]{Proposition}
\newtheorem{corollary}[theorem]{Corollary}
\theoremstyle{remark}
\theoremstyle{definition}
\newtheorem{remark}[theorem]{Remark}
\newtheorem{definition}[theorem]{Definition}

\newtheorem{problem}{Problem}
\newtheorem{conjecture}[problem]{Conjecture}

\title[Multidimensional continued fractions]{Constructing
multidimensional periodic continued fractions in the sense of
Klein}

\author{O.~N.~Karpenkov}

\date{7 April 2005}
\thanks{{\it AMS Subject Classification:} primary 11J70, secondary 11Y16.}
\thanks{Supported by SS-1972.2003.1 and RFBR-05-01-01012a grants.}

\keywords{Multidimensional continued fractions, convex polygons,
integer lattices}

\email[Oleg Karpenkov]{karpenk@mccme.ru}

\address{Poncelet Laboratory (UMI 2615 of CNRS and Independent University of Moscow)}

\begin{document}
\input epsf

\begin{abstract}{We consider the geometric generalization of ordinary continued
fraction to the multidimensional case introduced by F. Klein in
1895. A multidimensional periodic continued fraction is the union
of sails with some special group acting freely on these sails.
This group transposes the faces. In this article, we present a
method of constructing ``approximate'' fundamental domains of
algebraic multidimensional continued fractions and an algorithm
testing whether this domain is indeed fundamental or not. We give
some polynomial estimates on number of the operations for the
algorithm. In conclusion we present an example of fundamental
domains calculation for a two-dimensional series of
two-dimensional periodic continued fractions.}
\end{abstract}

\maketitle

\sloppy \normalsize

\tableofcontents

\section*{Introduction, definitions, and background}

The main goal of this paper is to introduce a new method for
constructing the fundamental domains of multidimensional periodic
continued fractions in the sense of Klein. The problem of
generalizing ordinary continued fractions to the
higher-dimensional case was posed by C.~Hermite~\cite{Herm} in
1839. A large number of attempts to solve this problem leads to
the birth of several different remarkable theories of
multidimensional continued fractions. In this paper we consider
the geometrical generalization of ordinary continued fractions to
the multidimensional case represented by F.~Klein in~1895 and
published by him in~\cite{Kle1} and~\cite{Kle2}.

A number of properties for ordinary continued fractions possesses
multidimensional analogies. H.~Tsuchihashi~\cite{Tsu} found the
connection between periodic multidimensional continued fractions
and multidimensional cusp singularities. J.-O.~Moussafir
in~\cite{Mou1} and O.~German in~\cite{Ger1} described
relationship between sails of multidimensional continued
fractions and Hilbert bases. M.~L.~Kontsevich and Yu.~M.~Suhov
discussed the statistical properties of the boundary of a random
continued fraction in~\cite{Kon}. The papers~\cite{Sku1}
and~\cite{Sku2} by B.~F.~Skubenko and~\cite{Ger2} by O.~N.~German
are dedicated to the generalization of ordinary continued
fractions with bounded above integer lengths of edges (the
numbers that corresponds to such continued fractions are the
numbers with the worst possible rational approximations). For the
classical theory of ordinary continued fractions we refer to the
book~\cite{Hin} by A.~Ya.~Hinchin. V.~I.~Arnold suggested to
investigate the geometry and combinatorics of continued fractions
(it contains the study of properties for such notions as affine
types of the faces of the sails, their quantities and
frequencies, integer angles between the faces, integer distances,
volumes and so on) in his article~\cite{Arn4} and the
book~\cite{Arn2}.

Some examples of the periodic continued fractions were calculated
in the papers~\cite{Kor0}, \cite{Kor2}, and~\cite{Kor3} by
E.~Korkina, \cite{Lac} and~\cite{Lac2} by G.~Lachaud, \cite{BP},
\cite{Par1}, \cite{Par1.1}, \cite{Par1.2}, and~\cite{Par2} by
A.~D.~Bruno and V.~I.~Parusnukov, \cite{Kar1} and~\cite{Kar2} by
the author. A nice collection of twodimensional continued
fractions was presented by K.~Briggs, see~\cite{site}.

\subsection{Definition of periodic multidimensional continued fractions}

In this section we recall some basic notions and definitions (see
also~\cite{Kar3}). Consider a space $\r^{n+1}$ ($n\ge 1$) over
$\r$. A point of $\r^{n+1}$ is said to be {\it integer} if all
its coordinates are integers. Two sets are called {\it
integer-affine $($integer-linearly$)$ equivalent} if there exists
an affine (linear) transformation of $\r^{n+1}$ preserving the
set of all integer points, and transforming the first set to the
second. A plane is called {\it integer} if it is integer-affine
equivalent to some plane passing through the origin and
containing the sublattice of the integer lattice, and the rank of
the sublattice is equivalent to the dimension of the plane. A
polyhedron is said to be {\it integer} if all its vertices are
integers.

Consider an integer plane and an integer point in the complement
to the plane. Let the Euclidean distance from the given point to
the given plane equal $l$. The minimal value of nonzero Euclidean
distances from integer points of the span of the the given plane
and the given point to the plane is denoted by $l_0$.  The ratio
$l/l_0$ is said to be the {\it integer distance} from the given
integer point to the given integer plane.

Consider arbitrary $n+1$ hyperplanes in~$\r^{n+1}$ that intersect
at a unique point, namely the origin. The complement to the union
of these hyperplanes consists of $2^{n+1}$ open orthants. Let us
choose an arbitrary orthant.
\begin{definition}
 The boundary of the convex hull of all integer
points except the origin in the closure of the orthant is called
the {\it sail} of the orthant. The set of all $2^{n+1}$ sails is
called the {\it $n$-dimensional continued fraction} constructed
according to the given $n+1$ hyperplanes.
\end{definition}

Two $n$-dimensional continued fractions are said to be {\it
equivalent} if the union of all sails of the first continued
fraction is integer-linear equivalent to the union of all sails
of the second continued fraction.

\begin{definition}
An operator in the group $SL(n+1,\z)$ is called an {\it integer
irreducible hyperbolic operator} if the following conditions hold:

i) the characteristic polynomial of this operator is irreducible
over $\q$;

ii) all its eigenvalues are distinct and real.
\end{definition}

Consider some integer irreducible hyperbolic operator $A\in
SL(n+1,\z)$. Let us take the $n$-dimensional spaces that span all
subsets of $n$ linearly independent eigenvectors of the operator
$A$. The spans of every $n$ eigenvectors uniquely define $n+1$
hyperplanes passing through the origin in general position. These
hyperplanes uniquely define the {\it multidimensional continued
fraction associated to $A$}.

\begin{definition}
An $n$-dimensional continued fraction associated to some integer
irreducible hyperbolic operator $A$ is called {\it an
$n$-dimensional continued fraction of an $(n+1)$-algebraic
irrationality}. The case of $n=1 (2)$ corresponds to {\it
one$($two$)$-dimensional continued fractions of quadratic
$($cubic$)$ irrationalities}.
\end{definition}

Now we formulate the notion of periodic continued fraction
associated an algebraic irrationality. Let $A$ be an integer
irreducible hyperbolic operator. Denote by $\Xi (A)$ the set of
all integer operators commuting with $A$. These operators form a
ring with standard matrix addition and multiplication. (As a
group $\Xi(A)$ is isomorphic to $\z^{n+1}$.)

Consider the subset of the set $SL(n+1,\z)\cap \Xi(A)$ that
consists of all operators with positive real eigenvalues and
denote it by $\bar\Xi (A)$. From the Dirichlet unit element
theorem (see.~\cite{BSh}) it follows that the subset $\bar
\Xi(A)$ forms a multiplicative Abelian group isomorphic to $\z^n$,
and that its action is free. Any operator of this group preserves
the integer lattice and the union of all $n+1$ hyperplanes, and
hence it takes the $n$-dimensional continued fraction onto itself
bijectively. (Whenever all eigenvalues are positive, the sails
are also taken onto themselves in a one-to-one way.) In addition,
the quotient of a sail under this group action is isomorphic to
an $n$-dimensional torus. These statements are based on the
generalization of the Lagrange theorem on ordinary continued
fractions. The combinatorial topological generalization of
Lagrange theorem was obtained by E.~I.~Korkina in~\cite{Kor1} and
its algebraic generalization by G.~Lachaud~\cite{Lac}.

\begin{theorem}{\bf (E.~I.~Korkina,~\cite{Kor1}.)}
Consider some orthant $C$ and the sail $K(C)$ corresponding to
this orthant. Suppose that there exists a combinatorial
isomorphism of the sail $V(C)$ that preserves the combinatorial
structure of the sail, then there exists an operator in the group
$GL(n+1,\z)$ taking the orthant $C$ and the sail $V(C)$ to
themselves and establishing the isomorphism.
\end{theorem}

Unfortunately, the proof of this theorem is not yet published.

The algebraic version of the generalized Lagrange theorem was
obtained by G.~Lachaud. The formulation of his theorem requires
special notation and definitions, so we just refer the reader to
the article~\cite{Lac}, since we will not use it further. For
more information on generalizations of the Lagrange theorem for
the ordinary continued fractions to the multidimensional case see
the papers~\cite{Kor1},~\cite{Tsu},~\cite{Lac} and~\cite{Lac3}.

By a {\it fundamental domain} of a sail we mean the union of some
faces that contains exactly one face from each equivalence class
(with respect to the action of the group $\bar\Xi (A)$).

\subsection{Different algorithms for constructing the sails
of multidimensional continued fractions}

A multidimensional periodic algebraic continued fraction is a set
of infinite polyhedral surfaces (i.e.,~sails), that contain an
infinite number of faces. As we have already mentioned, the
quotient of any sail under the Dirichlet group action is
isomorphic to an $n$-dimensional torus. The algebraic periodicity
of the polyhedron allows to reconstruct the whole continued
fraction knowing only the fundamental domain. Moreover, any
fundamental domain contains only a finite number of faces of the
whole algebraic periodic continued fraction. Hence we are faced
with the problem of finding a good algorithm that enumerates all
the faces for this domain.

There were no algorithm for constructing multidimensional
continued fractions until T.~Shintani's work~\cite{Shi} in~1976.
Let $F$ be a totally real algebraic field of degree $n$. We take
all different embeddings of $F$ into $\r$ and denote them by
$\varphi_i$, $i=1,\cdots,n$ (there are exactly $n$ different
embeddings, since $F$ is totally real). Consider the following
embedding of $F$ into $\r^n$. For an arbitrary element $x$ of $F$
we suppose
$$
x \rightarrow (\varphi_1(x),\varphi_2(x), \ldots, \varphi_n(x)).
$$
T.~Shintani considered the action of the group of all totally
positive elements for the ring of integers of $F$ (by
component-wise multiplication by totally positive integers $x_+$)
on $\r^n_+$ for the described embedding of $F$. He proved that
the fundamental domain for this action is the union of a finite
number of simplicial cones of special type. (Note that if we take
some other order for the embeddings $\varphi_{i'}$, then the
fundamental domains will be integer-linear equivalent to the
fundamental domains for the embeddings considered above.) The
statement of T.~Shintani and its proof is actually the basis for
the construction of one-dimensional continued fractions.
Following T.~Shintani's work, E.~Thomas and A.~T.~Vasques obtained
several fundamental domains for the two-dimensional case
in~\cite{Tho}. Finally, R.~Okazaki presented a method that permits
to construct fundamental domains for fields of arbitrary degree
in his article~\cite{Oka}. E.~Korkina
in~\cite{Kor0},~\cite{Kor2},~\cite{Kor3} and G.~Lachaud
in~\cite{Lac},~\cite{Lac2} produced an infinite number of
fundamental domains for periodic algebraic two-dimensional
continued fractions. The method used for constructing fundamental
domains of multidimensional continued fractions in these papers
was inductive. The method produces the fundamental domain face by
face, verifying that each new face does not lie in the same orbit
with some face constructed before. Applying the method, one can
find the fundamental domain in finitely many steps.

Later on J.~O.~Moussafir developed an essentially different
approach in his work~\cite{Mou2}. It works for an arbitrary (not
necessary periodic) continued fraction and computes any bounded
part of an infinite polyhedron. The approach is based on
deduction. One produces a conjecture on the face structure for a
big part of the continued fraction, then it remains to prove that
any conjectured face is indeed a face of the part. This method
can be also applied to the case of periodic continued fractions.

In the present paper we describe a new advanced deductive
construction adapted especially to fundamental domains of
periodic continued fractions. The construction involves a method
for conjecturing the structure of the fundamental domain and an
algorithm testing whether the conjectured domain is indeed
fundamental. The main advantage of our algorithm is the following:
the number of ''false'' vertices of our approximation is much
smaller than the number of ''false'' vertices of the
approximation in the method of J.~O.~Moussafir (so that the
computational time is considerable reduced).

Note that this algorithm substantially uses the periodicity of
the continued fraction and hence it is impossible to apply it to
non-periodic continued fractions.

We prove the following statement for the two-dimensional case.

{\it Suppose we have a conjecture on the structure of the
fundamental domain for some sail of a two-dimensional periodic
continued fraction. Let this domain contain $N$ faces of all
dimensions. The test of the conjecture $($our algorithm$)$
requires no more than $CN^4$ additions, multiplications and
comparisons of two integers, where $C$ is a universal constant
that does not depend on $N$.}

All previous verification algorithms work exponential time with
respect to $N$.

\begin{remark}
Here we do not take into account that the integers can be quite
large. To calculate upper bounds for the working time of the
algorithms we need to multiply upper bounds for the number of
operations by some polynomial of the coefficients of the matrix
defining the continued fraction.
\end{remark}

Using the present algorithm, the author both generalized almost
all known simple examples and series of examples of fundamental
domains constructed before, and found a lot of new examples and
series (see~\cite{Kar1} and~\cite{Kar2}). Using these examples,
the author found the complete list of all two-dimensional periodic
continued fractions constructed by matrices of small norm
($|*|\le 6$) up to the integer-linear equivalence relation,
see~\cite{Kar3}. By the norm of a matrix, here we mean the sum of
the absolute values of all its coefficients.

\subsection{Description of the paper}

This work is organized as follows. The new method of sail
construction consists of six steps. We discuss its plan in
Section~1. In section~2 we describe the two common steps for both
inductive and deductive methods. In this section we show how to
find the generators for the group of $SL(n,\z)$-matrices
commuting with the given one. All results of Section~2 are
well-known and are given for completeness of exposition (see also
the works~\cite{Coh} by Cohen and~\cite{Lac2} by G.~Lachaud). In
Sections~3 and~4, we discuss the essential new part of the method.
We show how to produce conjectures on fundamental domains in
Section~3. In Section~4 we describe the algorithm for conjecture
tests in the case of two-dimensional continued fractions. In that
section, we also say a few words about the higher dimensional
case. We conclude in Section~5 with the detailed study of one
example of the method's application (see also~\cite{Kar1}).

{\bf Acknowledgement.} The author is grateful to V.~I.~Arnold,
E.~I.~Korkina, G.~Lachaud, M.~A.~Tsfasman, and A.~B.~Sossinsky
for constant attention to this work and useful remarks, and to
the Institut de Math\'ematiques de Luminy (CNRS) for hospitality
and excellent working conditions.

\section{Description of the new construction}

\subsection{Outline of the new construction}
Now we briefly outline the main idea of the new construction of
one of the fundamental domains of the multidimensional continued
fraction corresponding to the given integer irreducible
hyperbolic operator. Suppose that we are given the integer
irreducible hyperbolic operator $A \in SL(n+1,\z)$. To compute
some fundamental domain of a sail of the continued fraction
associated to $A$ it is sufficient to do the following:

{\bf 1.} Compute a convex hull approximation of the sail. Namely
take a large enough convenient set of integer points and find its
convex hull.

{\bf 2.} Make a conjecture on some fundamental domain. Here we
need to guess a set of faces that might form a fundamental domain.
We do this by finding a repeatable pattern in faces geometry.

{\bf 3.} Prove the conjecture if possible.

{\bf 4.} If cannot prove the conjecture, start with {\bf 1.} but
with a larger convenient set of points. \vspace{2mm}

In this situation the following two questions are actual:\\
{\it How to find a convenient set of integer points for the
approximation
of the sail?}\\
{\it How to test whether the conjecture of a fundamental domain
of the sail is true or not?}\\
We give the answers to these questions in the present paper.

\subsection{Steps of the construction}

Let us briefly itemize the main steps of the method.

{\bf The deductive algorithm of constructing one of fundamental
domains for the sail of the given operator $A$ in the given
orthant.}

{\it Step 1.} Calculate the basis of the additive group of the
ring $\Xi (A)$.

{\it Step 2.} Calculate the basis of the group $\bar \Xi (A)$
(using the result of Step~1).

{\it Step 3.} Find some vertex of the sail.

{\it Step 4.} Make a conjecture on a fundamental domain of the
sail (using the results of Step~2 and Step~3).

{\it Step 5.} Test the produced (in Step~4) conjecture.

\begin{remark}
It is supposed that the fundamental domain conjectured in Step~5
and the basis $A_1, \ldots, A_n$ of the group $\bar\Xi(A)$
satisfy the following conditions:\\
i) the closure of the fundamental domain is homeomorphic to the disk;\\
ii) the operators $A_1, \ldots, A_n$ define the gluing of this
disk to the $n$-dimensional torus.
\end{remark}

Both inductive and deductive algorithms require the first and the
second steps. We describe these two steps in the next section.
All other steps are essential for our construction. In the method
by J.-O.~Moussafir~\cite{Mou2} the conjecture has been producing
using the approximation of the orthant by some rational orthant.
In the present paper we propose to produce conjectures for some
set of periods, see the description of Steps~3 and~4. We show how
to test conjectures in the case of two-dimensional continued
fractions in the description of Step~5. The result is partially
based on the theorem on integer-affine classification of
two-dimensional faces at the integer distances to the origin
greater than one from~\cite{Kar5}. (For the case of
$n$-dimensional continued fractions for $n\ge 3$, the last step
is quite complicated, since the classification of tree-dimensional
faces at the integer distances to the origin greater than one is
unknown.) In the last step we also investigate an important
particular case that seems to be quite common for periodic
$n$-dimensional continued fractions as well.

\begin{remark}
Note that all deductive algorithms are not algorithms in the
strict sense. One should choose some basis of $\bar \Xi (A)$ in
the right way, produce a good conjecture, and then test it. Even
the algorithmic recognition of the period for the given picture
of the boundary of the sail approximation is supposed to be a
hard problem. That is the reason why this ``algorithm'' cannot be
done completely by some computer program. But at the other hand,
the deductive algorithm is effective in practice. All of the
examples listed in the article~\cite{Kar1} were produced using
this algorithm. The examples of this paper generalize and expand
almost all known periods of the sails calculated before.
\end{remark}

\begin{remark}
The present method can be naturally generalized to the case of
Minkowski-Voronoi model of multidimensional continued fractions.
(For the definitions of Minkowski-Voronoi model see~\cite{Min}
and~\cite{Voro}.)
\end{remark}

\section{General questions concerning the lattice bases}

In this section we briefly discuss the questions which are
necessary for both inductive and deductive methods (Steps~1 and~2
of Section~1). The answers to these questions were known before
(see also~\cite{Lac2} and~\cite{Coh}).

\subsection{Step 1. Calculation of a basis of the additive group of the ring~$\Xi (A)$}

Let $V_0V_1 \ldots V_m$ be some tetrahedron with vertices $V_0$,
$V_1, \ldots, V_{m-1}$, and $V_m$. Denote by $P(V_0,V_1,\ldots,
V_m)$ the following parallelepiped:
$$
\Big\{V_0+\sum\limits_{k=1}^{m}\alpha _k \bar{V_0V_k}| 0\le
\alpha_k \le 1, k=1, \ldots m\Big\}.
$$

In this section we consider $\Xi(A)$ as an additive group. We
start the algorithm with the calculation of a basis for the group
$\Xi (A)$. Let us identify the space $Mat((n+1)\times (n+1) ,\r)$
with the space $\r^{(n+1)^2}$ and consider the standard metrics
for this space. So any integer operator corresponds to some
integer point, and the distance between two operators is the
Euclidean distance between the corresponding points in
$\r^{(n+1)^2}$. We consider a sum of absolute values of all
coefficients for some operator $A$ as a norm of the operator $A$
and denote it by $||A||$.

In Proposition~\ref{z^n} below we show that the set $\Xi (A)$ is
an additive group isomorphic to $\z^{n+1}$. Then by
Corollary~\ref{basis1} below it follows that there exists a basis
of the group $\Xi (A)$ contained in the parallelepiped
$P(0,E,A,A^2,\ldots,A^{n})$. Thus by Proposition~\ref{for_basis1}
below all norms of elements of the basis are bounded above by
$$
N'=\sum\limits_{i=0}^{n}||A^{i}||.
$$

\begin{remark}
Applying the LLL-algorithm described by A.~K.~Lenstra,
H.~W.~Lenstra and Jr., and L.~Lov\'asz in~\cite{LLL} to the
lattice generated by $0,E,A,A^2,\ldots,A^{n}$ one constructs some
reduced basis. This will decrease the number $N'$ to some
number~$N$. (It is not necessary to use LLL-algorithm here, put
just $N=N'$.)
\end{remark}

The set of integer operators contained in the parallelepiped is
also a subset of the following set:
$$
Mat((n+1)\times (n+1) ,\z)\cap B_N(O),
$$
where $B_N(O)$ --- is an $N$-neighborhood of the origin (i.e. the
ball of radius $N$ centered in the origin). The last set contains
less no more than $(2N+1)^{n+1}$ elements. This set contains all
integer operators of the parallelepiped
$P(0,E,A,A^2,\ldots,A^{n})$. Further we find the set of integer
points of the parallelepiped in a finite and polynomial with
respect to $N$ number of operations (i.e. we have to choose a
subset of all points of the parallelepiped in the set considered
above). Finally we choose some basis of the group $\Xi (A)$ using
the algorithm described in Proposition~\ref{bip}.

Now we formulate the statements mentioned above.

\begin{proposition}\label{z^n}
For any integer irreducible hyperbolic operator $A$ the set $\Xi
(A)$ forms an additive group isomorphic to $\z ^{n+1}$. \qed
\end{proposition}

The detailed proof of Proposition~\ref{z^n} can be found, for
instance, in the book~\cite{Lac2} by G.~Lachaud.

\begin{corollary}\label{basis1}
There exists a basis of the group $\Xi (A)$, such that all its
elements are contained in the parallelepiped
$P(0,E,A,A^2,\ldots,A^{n})$. \qed
\end{corollary}

\begin{proposition}\label{for_basis1}
Consider the parallelepiped $P(0,E,A,A^2,\ldots,A^{n})$. The
norms of all operators contained in this parallelepiped are
bounded above by
$$
\sum\limits_{i=0}^{n}||A^{i}||.
$$
\end{proposition}

The proof is straightforward. \qed

\begin{proposition}\label{bip}
Let there be given a maximal rank sublattice of the integer
lattice in the integer plane. Let $O;A_1, \ldots, A_n$ generate
this sublattice $($here $O$ is the origin of the lattice$)$. Then
there exists a basis $O;B_1, \ldots , B_n$ of the integer lattice
in the parallelepiped $P(O,A_1, \ldots ,A_n)$ such that for any
natural $i\le n$ the vertex $B_i$ belongs to the parallelepiped
$P(O,A_1,\ldots ,A_i)$.
\end{proposition}

\begin{proof}
We will inductively construct the basis $O;B_1, \ldots , B_n$. On
the $i$-th step we will construct the basis $O,B_1, \ldots , B_i$
inside the parallelepiped $P(O,A_1, \ldots ,A_i)$ that satisfies
all the conditions listed above for the lattice in the plane
spans the points $O,A_1, \ldots ,A_i$.

{\it Base of induction.} For $i=1$ we choose the closest to the
point $O$ integer point as $B_1$.

{\it Step of induction.} Suppose we have constructed vertices
$O;B_1,\ldots, B_{i-1}$ ($i\le n$) satisfying the induction
conditions. Now we construct an integer point $B_{i}$. Let the
integer distance between the point $A_{i}$ and the plane spanned
by the points $O,A_1, \ldots A_{i-1}$ be equal to $d_{i-1}$.
Consider the plane $\pi_i$ in the span of $O,A_1, \ldots ,A_{i}$
that is parallel to the plane spanned by $O,A_1, \ldots ,A_{i-1}$
and at the unit integer distance to that plane. The plane $\pi_i$
is integer. Therefore the intersection of $\pi_i$ with the
parallelepiped $P(O,A_1, \ldots, A_{i})$ contains at list one
integer point. We choose one of these points and denote it by
$B_{i}$.

As far as the spans of $O,A_1, \ldots ,A_{i-1}$ and $O,B_1,
\ldots ,B_{i-1}$ coincides, all integer points of the
parallelepiped $P(O,B_1, \ldots, B_{i})$ are at an integer
distance $0$ or $1$ to the plane spanning $O,B_1, \ldots
,B_{i-1}$. By induction assumption the points $O;B_1, \ldots
,B_{i-1}$ generates the integer sublattice of the corresponding
plane. Therefore the parallelepiped $P(O,B_1, \ldots ,B_{i-1})$
is empty (i.e. does not contain integer points different from
vertices of the parallelepiped). From the last two facts it
follows that the parallelepiped $P(O,B_1, \ldots ,B_{i})$ is
empty. Hence the points $O;B_1, \ldots ,B_{i}$ also generate the
integer sublattice of the corresponding plane.

So we have constructed the basis satisfying all the conditions of
the proposition by induction.
\end{proof}

The algorithm of this step does not seem to be the optimal one.
So the following question is natural here:

\begin{problem}
Find some effective algorithm of calculating an integer
sublattice for the $k$-dimensional plane of the space $\r ^m$ if
some basis $A_1, \ldots, A_k$ for some sublattice of this plane
is known $($note that we know nothing about the corresponding
quotient group$)$.
\end{problem}

\subsection{Step 2. Calculation of a basis of $\bar\Xi (A)$}

From the algorithmic point of view this step is the most
complicated. We describe only the idea for one of the simplest
algorithms here and give the corresponding references.

Let $\chi (x)$ be the characteristic polynomial of the operator
$A$ and let $\xi$ be one of the roots of $\chi (x)$. Consider the
following map
$$
h:\Xi(A)\rightarrow \q [\xi].
$$
For any element $B\in \Xi(A)$ there exists a unique representation
$B=p_B(A)$, where the degree of the polynomial $p_B$ is less than
$n+1$ (since the operators $E,A,\ldots A^{n}$ are linearly
independent). We put
$$
h (B)=p_B(\xi).
$$
Note that this map is an isomorphism between the ring $\Xi(A)$
and its image $h (\Xi(A))$. The addition and multiplication
operations in the image are induced by the addition and
multiplication operators of the field $\q [\xi]$
(see~\cite{Lac2}). Moreover the set $h (\Xi (A))$ forms an order
in the field $\q [\xi]$. By the Dirichlet unit theorem it follows
that there exist a basis for the units of this order and a number
$\rho$, such that the norms of all its elements are bounded above
by $\rho$. Since the method of constructing the constant $\rho$
is standard, we omit it. (For the construction of $\rho$ and
proofs we refer to~\cite{BSh}.) Note that the integer
$(n+1)$-dimensional volume of the symplex spanned by the basis
operators assign the minimal value. Now according to the
book~\cite{BSh} we construct the basis by enumeration of all
vectors of the set $h(\Xi(A))$ inside the ball $B_{\rho}(O)$,
where $B_{\rho}(O)$ is a $\rho$-neighborhood of the origin. The
preimage (i.e. $h ^{-1}$) of this basis gives us the basis of the
group of invertible elements in the ring $\Xi(A)$, and hence it
gives the basis of the subgroup $\bar\Xi (A)$.

\begin{remark}
The constant $\rho$ is extremely large (it equals the exponent of
some polynomial of the coefficients of the matrix $A$). The
effective algorithm for this step can be found in the
book~\cite{Coh} written by H.~Cohen. Using this algorithm one
finds the basis of units in the polynomial (with respect to the
coefficients of the matrix $A$) time.
\end{remark}

\begin{remark}\label{rem1}
Note that it is not necessary to find generators of $\bar\Xi (A)$.
The algorithm works for arbitrary $n+1$ linearly independent
operators of $\bar\Xi(A)$ (for more details see Remark~\ref{rem2}
below).
\end{remark}

\section{On fundamental domains and sail approximations}

The basis of the group $\bar \Xi (A)$ was calculated in the
previous section. Now we are coming to the main steps of the
algorithm. In this section we calculate one of the sail vertices
and show how to produce the conjectures.

\subsection{Step 3. How to calculate one of the sail vertices}

First let us find some integer point of the orthant containing
the sail. Consider an arbitrary orthant. Shift the standard unit
parallelepiped inside this orthant. Some integer point lies
inside the shifted parallelepiped. The coordinates of this point
coincide with integer parts of coordinates for one of $2^{n}$
vertices of this parallelepiped. Using this fact one can easily
find such point.

So we have found some integer point $P$ of the orthant, let us
find some vertex of the sail corresponding to this orthant.
Consider some integer plane $\pi$ passing through the origin such
that the intersection of $\pi$ with the orthant is at a unique
point (at the origin). Suppose that an integer distance from the
point $P$ to this plane is equal to $d$. Now we look through all
the symplexes obtained as intersections of our orthant with
parallel to $\pi$ planes at integer distances to the origin equal
$1, \ldots, d$. Suppose that the first symplex containing integer
points lies in the plane at an integer distance equal $d'\le d$.
The convex hull of all points of this symplex coincides with some
faces of the sail. All vertices of this face are vertices of the
sail. Choose an arbitrary one of them.

\subsection{Step 4. How to produce a conjecture on the fundamental domain of a sail}

Suppose that we know some point $V$ of the sail in the orthant,
and a basis $A_1, \ldots, A_n$ for the group $\bar\Xi(A)$. Now we
must produce a conjecture on a fundamental domain of the sail.
Let us briefly describe how to do this. First, we compute the set
of integer points that contains all vertices of some fundamental
domain of the sail. Secondly, we show how to choose the infinite
sequence of {\it special polyhedron approximation} for the sail.
Finally, using the picture of this approximations, we formulate a
conjecture on a fundamental domain of the sail.

\begin{proposition}
Let $V$ be a vertex of the sail of $n$-dimensional continued
fraction of a $(n+1)$-algebraic irrationality. Then there exists
fundamental domain of the sail such that all vertices of this
domain are contained in the convex hull $H$ of the origin and of
$2^n$ distinct points of the following form
$$
V_{\varepsilon_1,\ldots,\varepsilon_n}=\left(
\prod\limits_{i=1}^{n}A_i^{\varepsilon_i} \right) (V)
$$
where $\varepsilon_i\in \{0,1\}$ for $1\le i \le n$.
\end{proposition}

\begin{proof}
Consider the polyhedral cone $C$ with vertex at the origin and
base at the convex polyhedron with vertices
$V_{\varepsilon_1,\ldots,\varepsilon_n}$. We take the union of
all images of this polyhedral cone under the actions of operators
$$
A_{m_1,\ldots,m_n}=\prod\limits_{i=1}^{n}A_i^{m_i},
$$
for $1\le i \le n$, where $m_i\in \z$. Obviously this union is
equivalent to the union of the whole open orthant and the origin.
Therefore, any vertex of the sail is obtained from a vertex
contained in the cone $C$ by applying an operator
$A_{m_1,\ldots,m_n}$ for some integers $m_i$. Moreover, the
convex hull of all integer points of the given orthant contains
the convex hull of the vertices of the form
$A_{m_1,\ldots,m_n}(V)$. Hence the sail (i.e. the boundary of the
convex hull of integer points) is contained in the closure of the
complement in the orthant to the convex hull of all integer
points of the form $A_{m_1,\ldots,m_n}(V)$. This complement is a
subset of the union of polyhedra obtained from $H$ by an action
of some operator $A_{m_1,\ldots,m_n}$ (for some integers $m_i$,
$1\le i\le n$). This concludes the proof of the proposition.
\end{proof}

We skip the classical description of the computation of the
convex hull for the integer points contained in the polyhedron
$H$. Denote the vertices of this convex hull by $V_r$ for $0<r\le
N$ (here $N$ is the total number of such points).

\begin{definition}
The convex hull of the following finite set of points
$$
\{ A_{m_1,\ldots,m_n}(V_r)|1\le m_i \le m, \forall i :0 \le i \le
n \}
$$
is called the {\it $n$-th special polyhedron approximation} for
the sail.
\end{definition}

The defined set contains approximately and less than $m^nN$
points. (Since we calculated some image points for the boundary
of $H$ several times, we do not know the exact number of points.)
The number $N$ is fixed for the given generators $A_1, \ldots,
A_n$ and $m$ varies. We should try to make a good conjecture with
the least possible~$m$.

\begin{remark}
For all the examples listed in the paper~\cite{Kar1} (and for the
example of the last section) it was sufficiently to take $m=2$ to
produce the corresponding conjectures.
\end{remark}

\begin{remark}
Note that $N$ is a function defined on the set of all generators
of the group $\bar\Xi (A)$ and it does not depend on $m$.
Therefore the ``quality'' of the approximation also depends on
the choice of the basis.
\end{remark}

\begin{remark}\label{rem2}
Now we briefly discuss the case of Remark~\ref{rem1}. Suppose we
know operators $A_1, \ldots, A_n$ that generate only some full
rank subgroup of the group $\bar\Xi(A)$. Let the index of this
subgroup be equal to $k$. Then we are faced with the following
two problems. First the number $N$ will be approximately
$k$-times greater than in the previous case. Secondly one should
also find a conjecture on generators of the group $\bar\Xi(A)$.
\end{remark}

\section{Test of the produced conjectures in the two-dimensional case}

Now it remains to test the produced conjectures of Step~4. In
this section we explain how to test conjectures for the case of
two-dimensional periodic continued fractions. The test consists
of seven stages. It uses classification theorem from~\cite{Kar5}.
We prove here that these seven stages are sufficient for the
verification whether the produced conjecture is true or not. The
complexity of these stages is polynomial in the quantity all
faces.

\subsection{Brief description of the test stages and formulation
of the main results}

Suppose we  have a conjecture on some fundamental domain $D$ for
some sail of a two-dimensional periodic continued fraction
associated to some  integer irreducible hyperbolic operator $A$.
Also from Step~2 we know some basis $B_1, B_2$ of the group
$\bar\Xi(A)$. Let $p_k$ (for $k=0,1,2$) be the number of all
$k$-dimensional faces of the fundamental domain~$D$. Denote by
$F_i$ ($i=1,\ldots , p_2$) all two-dimensional faces, i.~e.
polygons. All adjacent to each face vertices and edges are known.
It is also conjectured that the fundamental domain $D$ and the
basis $B_1, B_2$
satisfy the following conditions:\\
i) the closure of the fundamental domain is homeomorphic to the two-dimensional disk;\\
ii) $B_1$ and $B_2$ define the gluing of this disk to the
$n$-dimensional torus (the fundamental domain $D$ is in
one-to-one correspondence with this torus).

{\it Test of the conjecture.} Our test consists of the following
seven stages:
\\
1. test of condition~i);
\\
2. test of condition~ii);
\\
3. calculation of all integer distances from the origin to the
two-dimensional planes containing faces $F_i$ and verification of
their positivity;
\\
4. test on nonexistence of integer points inside the pyramids
with vertices at the origin and bases at $F_i$ (here the integer
points of faces $F_i$ are permitted);
\\
5. test on convexity of dihedral angles (for all edges of the
fundamental domain);
\\
6. verification that all stars of the vertices are regular;
\\
7. test if all vertices of $D$ are in the same orthant.

\begin{theorem}\label{proverka}
The described conjecture test for the fundamental domain $D$
requires less than
$$
C(p_0+p_1+p_2)^4
$$
additions, multiplications and comparisons of two integers, where
$C$ is a universal constant that does not depend on~$p_i$.
\end{theorem}

\begin{remark}
Note that here we do not take into account the complexity of
additions, multiplications and comparison of two large integers.
We think of any such operation as of one operation  (as a unit of
time). There are known some linear with respect to the
coefficients of the matrix $A$ bounds for the number of digits of
such integers. So the complexity should be multiplied by some
polynomial of the coefficients of $A$.
\end{remark}

\begin{theorem}\label{dostatochno}
Let the set of faces $D$ satisfy the following conditions:
\\
1$)$ condition~i$)$;
\\
2$)$ condition~ii$)$;
\\
3$)$ positivity of all integer distances from the origin to the
two-dimensional planes containing faces $F_i$;
\\
4$)$ there are no integer points inside the pyramids with
vertices at the origin and bases at $F_i$ $($here the integer
points of faces $F_i$ are permitted$)$;
\\
5$)$ all dihedral angles are convex;
\\
6$)$ all stars of the vertices are regular;
\\
7$)$ all vertices of $D$ are contained in the same orthant.
\\
Then $D$ is a fundamental domain of some sail of the continued
fraction associated to the operator~$A$.
\end{theorem}

We start with the proof of Theorem~\ref{proverka}. Let us show
that all the stages listed above can be realized in the
polynomial time.

\subsection{Test of condition~i)}

First we need to test that the closures of any two
two-dimensional faces either do not intersect or intersect at a
vertex, or at an edge (and hence the closures of any two
one-dimensional faces either do not intersect or intersect at a
vertex). For this test we need to solve linear systems of two
equations and $2p_1$ inequalities that define faces $F_i$ (in
three variables). The number of such systems equals the number of
couples of faces, i.e. equals~$\frac{p_2(p_2-1)}{2}$. Since there
are only three variables and two linear equalities, any system
can be reduced to a system of inequalities in one variable in a
linear of $p_1$ time. This implies that we need no more than
$C_{1,1}p_1p_2^2$ operations to solve all systems, where
$C_{1,1}$ is some constant that does not depend on $p_i$.

Secondly we have to test that the edges are adjacent either to one
or to two (two-dimensional) faces. This can be done in less than
$p_1p_2$ single adjacency tests. Any adjacent test can be done in
a linear of $p_1$ time. This yields that we need no more than
$C_{1,2}p_1^2p_2$ operations for these tests, where $C_{1,2}$ is
some constant that does not depend on $p_i$.

Consider all (one-dimensional) edges of the closure of $D$ that
are adjacent to exactly one face of $D$. Thirdly we test that the
union of such edges is homeomorphic to the circle, i.~e. is
piece-wise connected and does not have self intersections. We
will check that this set is piece-wise connected, and that for
any vertex of the closure of $D$ either exactly two of the
described edges are adjacent to the vertex or none of the
described edges is adjacent to it. The test contains no more than
$p_1^2$ simple adjacency tests. Any simple adjacency test
requires a finite number of operations that does not depend on
$p_i$ for $i=0,1,2$. So, we need no more than $C_{1,3}p_1^2$
operations for these tests, where $C_{1,3}$ is some constant that
does not depend on $p_i$.

Further we will check that the union of the closures of all faces
of $D$ is piecewise connected. Let us find explicitly one of the
connected components. Consider an arbitrary face and all its
edges and vertices. Take all faces whose closures contain these
edges and vertices (except the first face) and consider all new
edges and vertices at their boundaries. We continue our
construction inductively. Suppose we have made $l$ steps and
constructed some part of the connected component. Consider again
all new (appeared on the $i$-th step) edges and vertices on its
boundary. Take all faces whose closures contain these edges and
vertices (except old faces). Since any edge of $D$ considers no
more than once, this algorithm requires no more than
$C_{1,4}'p_1$ simple adjacency tests. So we know one of the
connected components. Now it remains to check that all faces are
in this component and all edges and vertices are adjacent to
these faces. Finally we need no more than $C_{1,4}(p_0+p_1+p_2)$
operations.

Suppose all tests of this subsection are positive. So the closure
of $D$ is piecewise connected, and its boundary is homeomorphic
to the circle. The closure of $D$ is homeomorphic to the
two-dimensional disk iff its Euler characteristics equals one. We
need less than $4p_0+2p_1+p_2$ of additions for the calculation of
Euler characteristics for the closure of $D$.

We have proved the following lemma.
\begin{lemma}\label{pr1}
The first stage of the test requires no more than
$$
C_1(p^2_1p_2+p_1p_2^2+p_0)
$$
additions, multiplications and comparisons of two integers, where
$C_1$ is a universal constant that does not depend on~$p_i$. \qed
\end{lemma}

\subsection{Test of condition~ii)}

First we check that after the gluing of the disk $D$ we obtain a
nonsingular triangulated surface.

Let us test that there is no singular points inside the disk,
i.~e. open faces of the interior of the disk and their images
after gluing do not intersect. This requires $2(p_0+p_1+p_2)^2$
systems of linear equations and inequalities. Any such system can
be solved in a linear time with respect to $p_0$, $p_1$, $p_2$.
Hence we need no more than $C_{2,1}(p_0+p_1+p_2)^3$ operations,
where $C_{2,1}$ is a universal constant that does not depend
on~$p_i$.

Now we test that there is no singularities on edges on the
boundary of the disk. Here we need to check that for any edge of
the boundary there exists a unique edge that is gluing with the
first edge. This requires no more than $p_1^2$ adjacency tests in
$C_{2,2}$ operations each. Hence we need no more than
$C_{2,2}p_1^2$ operations ($C_{2,1}$ does not depend on~$p_i$).

Further we check that there is no singularities at the vertices
(after the gluing). For any vertex we need to check that the
union of all faces and edges containing the vertex in their
closure is homeomorphic to the disk. Now we know that faces of
the fundamental domain do not intersect. It remains only to test
that any edge is adjacent to exactly two faces (after gluing),
and that the union of all faces and edges containing the vertex
in their closure is piece wise connected and orientable. The test
of all these conditions requires no more than $hp_1p_2$ adjacency
tests for some constant $h$. Thus to test the vertex we need no
more than $C_{2,3}p_1p_2^2$ operations ($C_{2,3}$ does not depend
on~$p_i$). For all vertices of the boundary we need no more than
$C_{2,3}p_0p_1p_2^2$ operations.

Finally, let us check that this surface is homeomorphic to the
torus. First the Euler characteristic of the surface should be
zero (to verify this we need a linear of $p_i$ time). Secondly we
show the orientability of the surface. We orient the boundary
circle and check that any two boundary edge that we glue together
are glued with the opposite orientations. If for one of the
couples of such edges the orientation coincides, we have a Klein
bottle. Otherwise, we get a torus. This requires a linear of
$p_1$ number of operations. Therefore for this test we need less
than or equal $C_{2,4}(p_0+p_1+p_2)$ operations ($C_{2,4}$ does
not depend on~$p_i$).

It remains to show that the fundamental domain $D$ maps to the
obtained torus bijectively. For all faces this holds
automatically. The corresponding test for edges and vertices
requires $C_{2,5}(p_0+p_1)$ operations ($C_{2,5}$ does not depend
on~$p_i$).

We have proven the following lemma.
\begin{lemma}\label{pr2}
The second stage of the test requires no more than
$$
C_2(p_0p_1^2p_2+(p_0+p_1+p_2)^3)
$$
additions, multiplications and comparisons of two integers, where
$C_2$ is a universal constant that does not depend on~$p_i$. \qed
\end{lemma}

\subsection{Calculation of all integer distances from the origin
to the two-dimensional planes containing faces $F_i$, their
positivity}

Let $V_1(x_1,y_1,z_1)$, $V_2(x_2,y_2,z_2)$, and $V_3(x_3,y_3,z_3)$
be some integer points that do not lie in a straight line. Then
the following statement holds.

\begin{proposition}\label{distance}
An integer distance from the origin to the two-dimensional
integer plane that spans the points $V_1$, $V_2$, and $V_3$ is
equal to
$$
\frac{ \left| \left(
\begin{array}{ccc}
x_1 & x_2 & x_3 \\
y_1 & y_2 & y_3 \\
z_1 & z_2 & z_3 \\
\end{array}
\right) \right| } {|[\bar{V2-V1},\bar{V3-V1}]|_{\z}} ,
$$
where by $[\bar{V2-V1},\bar{V3-V1}]$ we denote the cross product
of the vectors $\bar{V2-V1}$ and $\bar{V3-V1}$ in $\r^3$, and by
${|[\bar{V2-V1},\bar{V3-V1}]|_{\z}}$ we denote the integer lengths
of the vector $[\bar{V2-V1},\bar{V3-V1}]$ $($i.~e. the greater
common divisor of the coordinates of this vector$)$. \qed
\end{proposition}

The proof is straightforward and is omitted here.

\begin{remark}
This proposition can be generalized to the higher dimensional
case.
\end{remark}

\begin{lemma}\label{pr3}
The third stage of the test requires no more than $C_3p_2$
additions, multiplications and comparisons of two integers, where
$C_3$ is a universal constant that does not depend on~$p_i$. \qed
\end{lemma}

\begin{proof}
Any integer distance to the origin can be found by the formula of
Proposition~\ref{distance} (this requires some finite number of
steps that does not depend on~$p_i$).
\end{proof}

\subsection{Test on nonexistence of integer points inside the pyramids
with vertices at the origin and bases at $F_i$}

First we formulate the following integer-linear classification
theorem.

\begin{theorem}\label{cor111}{\bf (See~\cite{Kar5}.)}
Any compact two-dimensional face of a sails of a two-dimensional
continued fraction contained in a plane at an integer distance to
the origin greater than one is integer-linear equivalent exactly
to one of the faces with vertices of the following list:

--- $(\xi,r-1,-r)$, $(a+\xi,r-1,-r)$, $(\xi,r,-r)$, where $a \ge 1$,
and $\xi$ and $r$ are relatively prime, and $r\ge 2$ and $0<\xi\le
r/2$;

--- $(2,1,b-1)$, $(2,2,-1)$, $(2,0,-1)$, where $b\ge 2$;

--- $(2,-2,1)$, $(2,-1,-1)$, $(2,1,2)$ and
$(3,0,2)$, $(3,1,1)$, $(3,2,3)$. \qed
\end{theorem}

For this stage we will also need the following proposition.

\begin{proposition}\label{tc_lin}
Consider two integer triangles in $\r^3$ with vertices
$A_1(x_{a_1},y_{a_1},z_{a_1})$, $A_2(x_{a_2},y_{a_2},z_{a_2})$,
$A_3(x_{a_3},y_{a_3},z_{a_3})$ and with vertices
$B_1(x_{b_1},y_{b_1},z_{b_1})$, $B_2(x_{b_2},y_{b_2},z_{b_2})$,
$B_3(x_{b_3},y_{b_3},z_{b_3})$. Let the plane of the first
triangle do not contain the origin. Then the triangle $A_1A_2A_3$
is integer-linear equivalent to the triangle $B_1B_2B_3$ $($in
the corresponding order$)$ iff the absolute value of the
determinant of the matrix $BA^{-1}$ equals one, and all the
coefficients of this matrix are integers, where
$$
A= \left(
\begin{array}{ccc}
x_{a_1} & x_{a_2} & x_{a_3}\\
y_{a_1} & y_{a_2} & y_{a_3}\\
z_{a_1} & z_{a_2} & z_{a_3}\\
\end{array}
\right) , \qquad B= \left(
\begin{array}{ccc}
x_{b_1} & x_{b_2} & x_{b_3}\\
y_{b_1} & y_{b_2} & y_{b_3}\\
z_{b_1} & z_{b_2} & z_{b_3}\\
\end{array}
\right) .
$$
\end{proposition}

We left the proof of this proposition as an easy exercise to the
reader.

\begin{lemma}\label{pr4}
The fourth stage of the test requires no more than $C_4p_2$
additions, multiplications and comparisons of two integers, where
$C_4$ is a universal constant that does not depend on~$p_i$.
\end{lemma}

\begin{proof}
Consider some face $F_i$. If an integer distance from the origin
to the plane containing the face equals one then all integer
points of the pyramid except the vertex at the origin are
contained in the base.

Suppose now that an integer distance from the origin to the plane
containing the face equals $r_i>1$. By Theorem~\ref{cor111} it
follows that this pyramid is triangular. First we calculate
integer lengths of the edges and the integer area of $F_i$ (in a
fixed number of operations). Further by Proposition~\ref{tc_lin}
the verification of an integer-linear type of the face is reduced
to the solution of nine integer linear equations on $\xi$ or to
the verification of some nine rationals to be integers.

Therefore for any face we need some constant number of operations
$C_4$ (where $C_4$ does not depend on~$p_i$). Thus the complexity
of the fourth stage of the test is no more than $C_4p_2$. This
concludes the proof of the lemma.
\end{proof}

\subsection{Test of the convexity of dihedral angles}
A dihedral angle is called {\it well-placed} if the origin is
contained in the corresponding opposite angle.

Note that the property of some angle to be well-placed is an
integer-linear invariant. Thus it is sufficient to check this
property for all dihedral angles for edges of the closure of $D$.
The test for each edge reduces to the solution of a system of
ordinary inequalities (without variables).

\begin{lemma}\label{pr5}
The fifth stage of the test requires no more than $C_5p_2^2$
additions, multiplications and comparisons of two integers, where
$C_5$ is a universal constant that does not depend on~$p_i$. \qed
\end{lemma}

\subsection{Verification that all 2-stars of the vertices are regular}

Let $p:W \to T^2$ be the universal covering of the torus after
the gluing. By the {$2$-star at the vertex of the universal
covering of the torus} we will call the union of all faces of the
universal covering (of dimensions no more than 2) to which the
given vertex is adjacent.

Let a vertex $v\in W$ map to $x\in \r^3$. The faces of the
universal covering to which the given vertex is adjacent maps to
the faces with the same property.

If $x\ne (a,0,0)$ for some positive $a$, then by $\bar{v}_n$ we
denote the vector $(1/n,0,0)$, $n \in \n$. (If $x=(a,0,0)$, then
by $\bar{v}_n$ we denote $(0,1/n,0)$.) A 2-star at $v$ is called
{\it regular} if for the sequence of rays $l_n$ passing through
the point $x+\bar{v}_n$ and having vertexes at the origin there
exists a positive $k$ such that for any $m\ge k$ the following
holds: the preimage ($p^{-1}$) of the intersection of the image
of 2-star at $v$ and $l_m$ in the universal covering consists of
exactly one point.

\begin{lemma}\label{pr6}
The sixth stage of the test requires no more than
$$C_6p_0p_1(p_0+p_2)$$
additions, multiplications and comparisons of two integers, where
$C_6$ is a universal constant that does not depend on~$p_i$.
\end{lemma}

\begin{proof}
For any face (edge) of the image of the 2-star for any vertex we
need to solve a system of no more than $p_1$ ($p_0$, respectively)
linear inequalities of the variable $\varepsilon =1/n$. Hence the
complexity of the test is no more than to
$$
C_6p_0(p_0p_1+p_1p_2),
$$
where $C_6$ is a universal constant that does not depend on~$p_i$.
\end{proof}

\subsection{Test for all vertices of $D$ to be in the same orthant;
conclusion of the proof of Theorem~\ref{proverka}}

\begin{lemma}\label{pr7}
The seventh stage of the test requires no more than $C_7p_0$
additions, multiplications and comparisons of two integers, where
$C_7$ is a universal constant that does not depend on~$p_i$.
\end{lemma}

\begin{proof}
We will use the following statement. Let $e$ be some nonzero
vector. The volume of the parallelepiped generated by the vectors
$e$, $A(e)$, and $A^2(e)$ equals zero iff the vectors $e$,
$A(e)$, and $A^2(e)$ generate some eigensubspace of $A$ (of
nonzero codimension).

Consider two vertices $x_1$ and $x_2$ of $D$. Let
$x(t)=tx_1+(1-t)x_2$. The vertices $x_1$ and $x_2$ are contained
in the same orthant iff the volume function $f(t)$ of the
parallelepiped generated by the vertices $x(t)$, $B_1(x(t))$, and
$B_1^2(x(t))$ does not have zeros in the segment $[0,1]$. Note
that $f(x)$ equals the determinant of the matrix generated by the
vectors $x(t)$, $B_1(x(t))$, and $B_1^2(x(t))$. Thus $f(x)$ is a
polynomial of the third degree with integer coefficients. All
solutions of $f(x)=0$ can be found explicitly. It remains to
compare them with 0 and 1.

Let us fix $x_1$ and vary $x_2$ (in the set of all vertices of
$D$). The corresponding test requires no more than $C_7p_0$
operations ($C_7$ does not depend on~$p_i$).
\end{proof}

\begin{remark}
As long as the integer operator $B_1$ is hyperbolic and
irreducible and all the coefficients of $f(x)$ are integers, it
follows that $f(x)$ has three distinct real roots. So the
question of existence of roots in $[0,1]$ can be reduced to
calculating critical points, comparing them with 0 and 1, and
comparing the critical values with $f(0)$ and $f(1)$. Here we
solve only one quadratic equation (instead of the cubic one
$f(x)=0$).
\end{remark}

{\it Conclusion of the proof of Theorem~\ref{proverka}.} From
Lemmas~\ref{pr1}---\ref{pr7} it follows that all seven stages of
the test require no more than $C(p_0+p_1+p_2)^4$ additions,
multiplications and comparisons of two integers, where $C$ is a
universal constant that does not depend on~$p_i$. The proof of
Theorem~\ref{proverka} is complete. \qed

\begin{remark}
Actually a stronger statement holds. The whole test requires no
more than $\bar{C}(p_0+p_1+p_2)^3$ additions, multiplications and
comparisons of two integers, and $\tilde{C}(p_0+p_1+p_2)^4$
logical operations, where $\bar{C}$ and $\tilde{C}$ are universal
constants that do not depend on~$p_i$.
\end{remark}

It remains to prove that this seven stages are sufficient for the
test.

\subsection{Lemma on the injectivity of the face projection}
We prove Theorem~\ref{dostatochno} in four lemmas.

First let us give the necessary notation. Let the operators $B_1$
and $B_2$ generate $\bar \Xi (A)$. For any integers $n$, $m$ by
$B_{n,m}$ we denote an operator $B^n_1B^m_2$. We suppose that our
domain $D$ satisfies conditions~1---7 of
Theorem~\ref{dostatochno}. Let
$$
U=\bigcup\limits_{n,m\in \z} B_{n,m}(D).
$$

Consider the unit two-dimensional sphere $S^2$ centered at the
origin $O$. We denote by $\pi$ the following map:
$$
\pi : \r^3 \setminus {O} \to S,
$$
where any point $x \in \r^3 \setminus {O}$ maps to the point at
the intersection of $S^2$ and the ray with vertex at the origin
and containing $x$.

\begin{lemma}\label{proverka1}
For any face of the polygonal surface $U$ the map $\pi$ is
well-defined and injective on it.
\end{lemma}

\begin{proof}
Consider any two-dimensional face $F$ of the surface $U$. By
condition~3, the distance from the origin to the plane containing
$F$ is greater than zero. Hence this plane does not contain the
origin. Then $\pi$ is well-defined and injective on $F$.

Let now $E$ be some edge of $U$. By conditions~1 and~2, this edge
is adjacent to some two-dimensional face and therefore is
contained in some plane  that does not pass through the origin.
Thus the line containing $E$ does not pass through the origin.
Hence $\pi$ is well defined and injective on $E$.

The injectivity for the vertices is obvious.
\end{proof}

\subsection{Lemma on the finite covering of the fundamental domain}

Let $x \in \r^3 \setminus {O}$. Denote by $N_x$ the tetrahedral
angle with vertex at the origin and base with vertices $x$,
$B_1(x)$, $B_1B_2(x)$, and $B_2(x)$. Notice that
$$
\left( \bigcup\limits_{n,m\in \z} B_{n,m}(N_x) \right) \setminus O
$$
is one of eight orthants of the continued fraction associated to
$A$ that contains $x$. Note that from conditions~1, 2, and~6 it
follows that all points of $D$ are contained in one open orthant,
we denote it by $K$.

\begin{lemma}\label{proverka2}
Let $x$ be some point of the open orthant $K$. Then the union of
all faces of $D$ is contained in a finite union of solid angles
of the type $B_{n,m}(N_x)$.
\end{lemma}

\begin{proof}
By the Dirichlet unit theorem~\cite{BSh} it follows that for any
interior point $a$ of the open orthant $K$ there exists an open
neighborhood satisfying the following condition. The neighborhood
can be covered by four solid angles of the type $B_{n,m}(N_x)$
when $a$ belongs to an edge of some $B_{k,l}(N_x)$; by two solid
angles when $a$ belongs to the face of some $B_{k,l}(N_x)$; and
by one solid angle in the remaining cases. In any case the
neighborhood can be covered by some finite union of solid angles
of the type $B_{n,m}(N_x)$.

Consider a covering of $D$ by such neighborhoods that correspond
to each point of the closure of $D$. Since the closure of $D$ is
closed and bounded in $\r ^3$, it is compact. Hence this covering
contains some finite subcovering. Therefore the union of all
faces of $D$ is contained in the finite union of solid angles of
type $B_{n,m}(N_x)$. The proof is complete.
\end{proof}

\begin{corollary}\label{proverka3}
Let $x$ be contained in the open orthant $K$. Then the solid
angle $N_x$ contains only points from a finite number of
fundamental domains of the type $B_{n,m}(D)$.
\end{corollary}

\begin{proof}
From the last lemma it follows that $D$ is contained in the finite
union $\bigcup\limits_{k=1}^{l}B_{n_k,m_k}(N_x)$ (for some
positive $l$). Then the solid angle $N_x$ can contain only points
of the fundamental domains $B_{-n_k,-m_k}(D)$ for $1 \le k \le l$.
\end{proof}

\subsection{Lemma on the bijectivity of the projection}

\begin{lemma}\label{proverka4}
The map $\pi$ bijectively takes the polygonal surface $U$ to the
set $S^2\cap K$.
\end{lemma}

\begin{proof}
As it was shown above, the surface $U$ is contained in $K$ and is
taken to $S^2\cap K$ under the map $\pi$.

Let us introduce the following notation. By condition~2, the
operators $B_1$ and $B_2$ glue the fundamental domain into the
torus $T^2$. Let $W$ be the universal covering of $T^2$. The face
decomposition on $T^2$ lifts to a face decomposition on $W$.
There is a natural two-parametric family (with two integer
parameters) of projections $p_{n,m}:W \to U$ that maps faces to
faces (since the group of shifts $B_{k,l}$ acts on $U$). Let us
choose one of these projections and denote it by $p$ ($p:W \to
U$).

Consider the map $\pi\circ p : W \to S^2 $. This map does not
have branch points at the images of open faces of $W$, since any
face of $W$ bijectively maps to some face of $U$, and the
corresponding face of $U$ injectively maps to $S^2 \cap K$ by
Lemma~\ref{proverka1}.

Two faces with common edge of $W$ map to some two faces with
common edge of $U$, such faces of $U$ generate a well-placed
dihedral angle, and hence also injectively map to $S^2 \cap K$.
So the map $\pi \circ p$ does not have branch points at the
images of open edges.

Now we consider some vertex $v$ of $W$. The edges and faces of
$W$ with common vertex $v$ by condition~6 form a regular $2$-star.
This edges also maps to some edges of $U$ with common vertex.
Thus there exist a sequence of points that tends to $\pi \circ p
(v)$  (contained in $S^2$) such that the preimage of any point of
the sequence has exactly one preimage in the $2$-star of $v$.
Hence $\pi \circ p$ does not have branch points at the sheet
containing the star at $\pi \circ p (v)$. Therefore, $\pi \circ p$
does not have any branch points at the vertices.

So the map $\pi \circ p: W \to S^2 \cap K$ does not have branch
points.

Consider an arbitrary point $x \in S^2 \cap K$ and the solid
angle $N_x$ corresponding to it. Let $x_1$ and $x_2$ be two
points of $S^2 \cap N_x$. Now we will show that the preimages
$(\pi \circ p)^{-1}(x_1)$ and $(\pi \circ p)^{-1}(x_2)$ contain
the same number of points. Let us join the points $x_1$ and $x_2$
by some curve inside $S^2 \cap N_x$. By
Corollary~\ref{proverka3}, we know that the preimage of this
curve is contained in a finite number of faces of $W$. Since
there are no branch points in any face (and their number is
finite) and there are no boundary faces of $W$, the number of
preimages for $\pi \circ p$ is some (finite) discrete and
continuous function on this curve. Therefore the number of
preimages for $\pi \circ p$ of any two points of $S^2 \cap N_x$
is the same. Hence the number of preimages for $\pi \circ p$ of
any two points of $S^2 \cap K$ is the same.

From this we conclude that the projection $\pi \circ p$ of the
universal covering $W$ (homeomorphic to an open disk) to $S^2
\cap K$ (i.~e. homeomorphic to an open disk) is a nonramified
covering with finitely many sheets. Since the covering (of an
open disk by an open disk) is piece-wise connected, the number of
sheets equals one. Since by definition $p: W \to U$ is surjective
and by the all above it is injective, the maps $p: W \to U$ and
$\pi :U \to S^2 \cap K$ are bijective. This concludes the proof
of the lemma.
\end{proof}

\subsection{Lemma on convexity}

Since any ball centered at the origin contains only a finite
number of vertices of $U$ (and they do not form a sequence
tending to the boundary of the orthant $K$), the polyhedral
surface $U$ divide the space $\r^3$ into two connected components.
Denote by $H$ the connected component of the complement to $U$
that does not contain the origin.

\begin{lemma}\label{proverka5}
The set $H$ is convex.
\end{lemma}

\begin{proof}
Suppose that some plane passing through the origin intersects the
polygonal surface $U$ and does not contain any vertex of $U$. By
Lemma~\ref{proverka4} such plane intersects $U$ at some piecewise
connected broken line with an infinite number of edges. The
complement of the plane to this broken line consists of two
connected components. By assumption all vertices of this broken
line are contained in open edges of $U$. By condition~5 all
dihedral angles of $U$ are well-placed. Thus the angle at any
vertex of intersection of $H$ with our plane is less than
straight angle. Hence by the previous lemma the intersection is
convex.

Consider the set of all planes that pass through the origin,
intersect $U$ and do not contain vertices of $U$. This set is
dense in the set of all planes passing through the origin and
intersecting $U$. Therefore by the continuity reasons it follows
that the intersection of $H$ with any plane passing through the
origin (and intersecting $U$) is convex.

Now we prove that the set $H$ is convex. Let  $x_1$ and $x_2$ be
some points of $H$. Consider the plane that spans $x_1$, $x_2$,
and the origin. This plane intersects $U$ since $x_1$ is in $H$
and the origin is not in $H$. By the above the intersection of
$H$ with this plane is convex. Hence the segment with endpoints
$x_1$ and $x_2$ is contained in $H$. Thus $H$ is convex (by the
definition of convexity).

Lemma~\ref{proverka5} is proven.
\end{proof}

\subsection{Conclusion of the proof of Theorem~\ref{dostatochno}:
the main part}

So the constructed polygonal surface $U$ possess the following
properties:
\\
--- by Lemma~\ref{proverka5} $U$ bounds the convex set $H$;
\\
--- by construction all vertices of $U$ are integer points;
\\
--- by condition~4 the set $K\setminus H$ does not contain integer points.

Therefore the polygonal surface $U$ is the boundary of the convex
hull of all integer points inside $K$. Thus by definition $U$ is
one of the sails of the continued fraction associated to the
operator~$A$. This concludes the proof of
Theorem~\ref{dostatochno}. \qed

Let us formulate one important conjecture here.

\begin{conjecture}
Conditions~1---6 imply condition~7.
\end{conjecture}

\subsection{On the verification of the conjecture for the multidimensional case}

Here we briefly outline an idea how to test the conjecture for
fundamental domains of multidimensional continued fractions.

{\it Conjecture for the multidimensional case.} Suppose we have a
conjecture on some fundamental domain $D$, and also some basis
$B_1, \ldots ,B_n$ of the group $\bar\Xi(A)$. Let also the
fundamental domain and the basis posses the
following properties:\\
i) the closure of the fundamental domain is homeomorphic to the disk;\\
ii) the operators $B_1, \ldots B_n$ define the gluing of this disk
to the $n$-dimensional torus.

{\it How to test the conjecture for fundamental domains of
multidimensional continued fractions?} The verification of
conditions i) and ii) is straightforward and is omitted here. If
these conditions hold we verify if all the $n$-dimensional faces
of the fundamental domain are faces of the sail. It can be done
in the following way.

Suppose that integer distances from the origin to the planes of
faces $F_i$ are equal to $d_i$ ($i=1,\ldots , p$, where $p$ ---
is the number of all $n$-dimensional faces). Our conjecture is
true iff for all $i=1, \ldots , p$ the following
conditions hold.\\
a) For any integer $d<d_i$ consider the plane parallel to the
face $F_i$ with integer distances to the origin equals $d$. The
intersection of our orthant with this plane
does not contain any integer point.\\
b) For $d=d_i$ the convex hull of all integer points
in the intersection coincides with face $F_i$.\\
The verification of conditions~a) and~b) is quite complicated
from the algorithmical point of view.

We conclude this section with the important inverse question of
construction periodic continued fractions.

\begin{problem}{\bf (V.~I.~Arnold.)}
Does there exist an algorithm to decide whether a given type of
fundamental domain is realizable by a periodic continued fraction.
\end{problem}

The answer to this question is unknown even for the
two-dimensional periodic continued fractions.

\section{An example of the calculation of a fundamental domain}

The example of two-dimensional sails in this section was announced
in the article~\cite{Kar1} by the author. For arbitrary integer
numbers $m$ and $n$ we denote by $A_{m,n}$ the Sylvester operator
$$
\left(
\begin{array}{ccc}
0 &1 &0 \\
0 &0 &1 \\
1 &-m &-n \\
\end{array}
\right) .
$$
We construct fundamental domains for some particular
two-dimensional subfamily of the two-dimensional continued
fractions corresponding to Sylvester operators for the orthant
that contains the point $(0,0,1)$.

\begin{theorem}\label{t31}
Let $m=b-a-1$, $n=(a+2)(b+1)$ $($where $a,b\ge 0)$. Consider the
sail of the operator $A_{m,n}$ containing the point $(0,0,1)$.
Let $A=(1,0,a+2)$, $B=(0,0,1)$, $C=(b-a-1,1,0)$, and
$D=((b+1)^2,b+1,1)$.
Then the following set of faces forms one of the fundamental domains:\\
1) the vertex $A$;\\
2) the edges $AB$, $AD$, and $BD$;\\
3) the triangular faces $ABD$ and $BDC$.

\end{theorem}

\begin{figure}[h]
$$\epsfbox{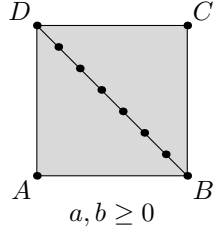}$$
\caption{The closure of the fundamental domain of a sail of a
fraction associated to the operator $A_{b-a-1,(a+2)(b+1)}$ (here
$b=6$, and $a$ is arbitrary).}\label{fig2.3}
\end{figure}

The closure of the fundamental domain is homeomorphic to the
square shown on Fig.~\ref{fig2.3} (for the case of an arbitrary
$a$, and $b=6$).

{\it Proof.} {\it Steps~1 and~2.} We omit the first and the
second steps (these steps are classical, see~\cite{Coh})) and
here write down the result. The following two operators generate
the group $\bar\Xi (A)$:
$$
X_{a,b}=A_{m,n}^{-2}, \quad
Y_{a,b}=A_{m,n}^{-1}\big(A_{m,n}^{-1}-(b+1)I\big),
$$
where $I$ is the identity element in the group $SL(3,\z)$.

{\it Step~3.} We prove that $(0,0,1)$ is a vertex of the sail.
Consider the plane passing through $A$, $B$, and $D$:
$$
(-1-a)x+(ab+a+b+1)y+z=1.
$$
As far as the equations (in variables $x$, $y$, and $z$)
$$
(-1-a)x+(ab+a+b+1)y+z=\alpha
$$
do not have any integer solution for $0<\alpha < 1$, an integer
distance from $ABD$ to the origin is equal to one. There are
exactly three integer points ($A$, $B$, and $D$) in the
intersection of the plane and the orthant (we left the proof of
that fact for the reader as an exercise).

{\it Step 4.} The conjecture of the fundamental domain was
produced in the statement of this theorem.

So it remains to complete {\it Step~5}: to test the conjectured
fundamental domain. For the test we need some extra points:
$$
\begin{array}{l}
E=X_{a,b}^{-1}(B)=(1,-ab-a-2b-2, \\
\makebox[2.5cm]{}a^2b^2+2a^2b+4ab^2+a^2+8ab+4b^2+5a+7b+5);\\
F=Y_{a,b}(B)=(-a-2,1,0);\\
H=X_{a,b}^{-1}(F)=(0,-b-1,ab^2+2ab+2b^2+a+4b+3).\\
\end{array}
$$

{\it 1.} ({\it{}Test of condition} i).) It can be shown in the
usual way that the faces have the common edge $BD$, and the edges
intersect only at vertices. This implies that all adjacencies are
correct, and that only one or two faces are adjacent to each edge.
The closure of the boundary is a closed broken line $ABCDA$,
homeomorphic to the circle.

{\it 2.} ({\it{}Test of condition} ii).) By direct calculations
it follows that by the operator $X_{a,b}$ action the segment $AB$
is taken to the segment $DC$ (the point $A$ maps to the point $D$
and $B$ to $C$) and by the operator $Y_{a,b}$ action the segment
$AD$ is taken to the segment $BC$ (the point $A$ maps to the
point $B$ and $D$ to $C$). Obviously that no other points glue
together. The Euler characteristic of the obtained surface equals
$2-3+1$, i.e. zero, and the surface is orientable.

{\it 3.} ({\it{}Calculation of all integer distances from the
origin to the two-dimensional planes containing faces}.) Let us
calculate integer distances from the origin to the two-dimensional
planes containing faces $ABD$ and $BDC$ by the formula of
Lemma~\ref{distance}. An integer distance to the plane of $ABD$
equals
$$
\frac{1}{b+1} \cdot \left| \left(
\begin{array}{ccc}
1 & 0 & b^2+2b+1\\
0 & 0 & b+1\\
a+1 & 1 & 1\\
\end{array}
\right) \right| =\frac{b+1}{b+1} =1.
$$

An integer distance to the plane of $BDC$ equals
$$
\frac{1}{b+1} \cdot \left| \left(
\begin{array}{ccc}
0 & b^2+2b+1 & b-a-1\\
0 & b+1 & 1\\
1 & 1 & 0\\
\end{array}
\right) \right| =\frac{ab+2b+a+2}{b+1} =a+2.
$$

{\it 4.} ({\it{}Test on nonexistence of integer points inside the
pyramids with vertices at the origin and bases at the faces}.)
Since the integer distance from the origin to the plane
containing $ABD$ equals one, the pyramid corresponding to $ABD$
does not contain integer points different from $O$ and the points
of the face $ABD$.

The face $BDC$ is integer-linear equivalent to the face with
vertices $(1,a+1,-a-2)$, $(b+2,a+1,-a-2)$, $(1,a+2,-a-2)$ of the
list of Theorem~\ref{cor111}. The corresponding transformation
taking $BCD$ to the face of the list of Theorem~\ref{cor111} is
the following:
$$
\left(
\begin{array}{ccc}
b+1 & b-a-1 & b-a\\
1 & 1 & 1\\
0 & -1 & -1\\
\end{array}
\right) .
$$
By Theorem~\ref{cor111} the pyramid corresponding to $BDC$ does
not contain integer points different from $O$ and the points of
the face $BDC$.

{\it 5.} ({\it{}Test on convexity of dihedral angles}.) Let us
first consider the edge $BD$. This edge is adjacent to the faces
$ABD$ and $BDC$. The face $ABD$ is contained in the plane
$f_{ABD}(x,y,z)=0$, and the face $BDC$ is contained in the plane
$f_{BDC}(x,y,z)=0$, where
$$
\begin{array}{ccl}
f_{ABD}(x,y,z) & = & (-1-a)x+(ab+a+b+1)y+z-1;\\
f_{BDC}(x,y,z) & = & x+(b+1)y-(a+2)z+(a+2).\\
\end{array}
$$
To test that the dihedral angle corresponding to the edge $BD$ is
well-placed it is sufficient to verify the following: the point
$C$ and the origin $O$ lie in different half-spaces with respect
to the plane that spans the points $A$, $B$, and $D$; the points
$A$ and $O$ lie in different half-spaces with respect to the plane
that spans the points $C$, $B$, and $D$. So we need to solve the
following system:
$$
\left\{
\begin{array}{ccc}
f_{ABD}(C)\cdot f_{ABD}(O) & < & 0\\
f_{BDC}(A)\cdot f_{BDC}(O) & < & 0\\
\end{array}
\right. .
$$
This system is equivalent to the following one:
$$
\left\{
\begin{array}{rcc}
(a^2+3a+2)\cdot (-1) & < & 0\\
(-a^2-3a-1)\cdot (a+2) & < & 0\\
\end{array}
\right. .
$$

Since $a\ge 0$, the inequalities hold. Thus the dihedral angle
associated with the edge $BD$ is well-placed.

Since the cases of dihedral angles associated to the edges $AB$
(and the faces $ADB$ and $AEB$) and $BC$ (and the faces $BDC$ and
$CBF$) can be verified in the same way, we omit their
descriptions.

This concludes the test of condition~5.

{\it 6.} ({\it{}Verification that all 2-stars of the vertices are
regular}.) There is only one vertex in the torus decomposition.
Any lift of this point to the universal covering $W$ is adjacent
to six edges and six faces. Consider a vertex of the universal
covering that maps to the point $B$. The corresponding $2$-star
maps to six edges $BC$, $BD$, $BA$, $BE$, $BH$, and $BF$ and to
six faces $BCD$, $BDA$, $BAE$, $BEH$, $BHF$, and $BFC$, where
$$
H=X_{a,b}^{-1}(F)=(0,-b-1,ab^2+2ab+2b^2+a+4b+3).
$$
We will check that  for any sufficiently small positive
$\varepsilon$ a ray $l_\varepsilon$ with vertex at the origin and
passing through the point $P_\varepsilon=(\varepsilon, 0, 1)$
intersects the exactly one of the faces of $2$-stars.

First we will check that for any sufficiently small positive
$\varepsilon$ the ray $l_\varepsilon$ intersects the triangle
$BCF$. Or, equivalently, that the ray $l_\varepsilon$ is
contained in the trihedral angle with vertex at the origin $O$ and
base in the triangle $BCF$. The two-dimensional face of the
trihedral angle containing $B$, $C$, and $O$ can be defined by
$f_{ABO}=0$; the two-dimensional face of the trihedral angle
containing $B$, $F$, and $O$ can be defined by $f_{BFO}=0$; the
two-dimensional face of the trihedral angle containing $C$, $F$,
and $O$ can be defined by $f_{CFO}=0$, where
$$
\begin{array}{ccl}
f_{BCO}(x,y,z) & = & x+(a+1-b)y;\\
f_{BFO}(x,y,z) & = & x+(a+2)y;\\
f_{CFO}(x,y,z) & = & z.\\
\end{array}
$$
For any sufficiently small positive $\varepsilon$ the ray
$l_\varepsilon$ is contained in the dihedral angle defined above
if the following conditions hold: the points $P_\varepsilon$ and
$F$ are in the same closed half-space with respect to the plane
$f_{BCO}=0$; the points $P_\varepsilon$ and $C$ are in the same
closed half-space with respect to the plane  $f_{BFO}=0$; the
points $P_\varepsilon$ and $B$ are in the same closed half-space
with respect to the plane  $f_{CFO}=0$. Since the points
$P_\varepsilon$ and $B$ are close to each other for sufficiently
small $\varepsilon$, they are in the same closed half-space with
respect to the plane  $f_{CFO}=0$. Now we check the remaining two
conditions:
$$
\left\{
\begin{array}{ccc}
f_{BCO}(P_\varepsilon)\cdot f_{BCO}(F) & \ge & 0\\
f_{BFO}(P_\varepsilon)\cdot f_{BFO}(C) & \ge & 0\\
\end{array}
\right. \quad \Leftrightarrow \quad \left\{
\begin{array}{rcc}
(-b-1)\varepsilon & \ge & 0\\
(b+1)\varepsilon & \ge & 0\\
\end{array}
\right. .
$$
Since $b, \varepsilon \ge 0$ the first inequality does not hold.
Thus for any sufficiently small positive $\varepsilon$ the ray
$l_\varepsilon$ does not intersect the triangle $BCF$.

The cases of the triangles $BCD$, $BDA$, $BAE$, $BEH$, and $BHF$
are similar to those described above and are omitted here.

The ray $l_\varepsilon$ (for any sufficiently small positive
$\varepsilon$) intersects the bijective image of a $2$-star of
the vertex at exactly one point contained in the edge $AB$.
Therefore all $2$-stars associated to the vertices are regular.

{\it 7.} ({\it Test that all the vertices of $D$ are in the same
orthant}.) The test of the seventh stage for this theorem is
trivial since $D$ contains exactly one vertex.

This concludes the proof of Theorem~\ref{t31}. \qed


\begin{thebibliography}{99}
\bibitem{Arn2}
V.~I.~Arnold, {\it Continued fractions}, M.: Moscow Center of
Continuous Mathematical Education, (2002).
\bibitem{Arn4}
V.~I.~Arnold, {\it Higher dimensional continued fractions},
Regular and Chaotic Dynamics, v.~3(3), pp.~10--17, (1998).
\bibitem{BSh}
Z.~I.~Borevich, I.~R.~Shafarevich, {\it Number theory}, 3 ed, M.,
(1985).
\bibitem{site}
K.~Briggs, {\it Klein polyhedra},
http://keithbriggs.info/klein-polyhedra.html, (2002).
\bibitem{BP}
A.~D.~Bryuno, V.~I.~Parusnikov, {\it Klein polyhedrals for two
cubic Davenport forms}, Mathematical notes, 56(4), (1994), 9-27.
\bibitem{Coh}
H.~Cohen, {\it A Course in Computational Algebraic Number Theory},
Graduate texts in mathematics. Berlin, Springer, (1973).
\bibitem{Ger1}
O.~N.~German, {\it Sails and Hilbert Bases}, Proc. of Steklov
Ins. Math, v. 239(2002), pp. 88-95.
\bibitem{Ger2}
O.~N.~German, {\it Sails and Norm Minima of Lattices}, To appear
in Matem. sbornik (2004).
\bibitem{Herm}
C.~Hermite, {\it Letter to C.~D.~J.~Jacobi}, J. Reine Angew.
Math. vol.~40, (1839), p.~286.
\bibitem{Hin}
A.~Ya.~Hinchin, {\it Continued fractions}, M.: FISMATGIS, (1961).
\bibitem{Kar1}
O.~Karpenkov, {\it On tori decompositions associated with
two-dimensional continued fractions of cubic irrationalities},
Func. an. and appl., v.~38(2004), no~2, pp.~28-37.
\bibitem{Kar2}
O.~N.~Karpenkov, {\it On two-dimensional continued fractions for
integer hyperbolic matrices with small norm,} Uspehi Mat. Nauk,
vol.~59(2004), no.~5, pp.~149--150.
\bibitem{Kar3}
O.~N.~Karpenkov, {\it On examples of two-dimensional periodic
continued fractions}, preprint, Cahiers du Ceremade, UMR 7534,
Universit\'e Paris-Dauphine, (2004).
\bibitem{Kar5}
O.~N.~Karpenkov, {\it Classification of three-dimensional
multistory completely empty convex marked pyramids}, Uspehi Mat.
Nauk, vol.~60, (2005), no.~1, pp.~169--170.
\bibitem{Kle1}
F.~Klein, {\it Ueber einegeometrische Auffassung der gew\"ohnliche
Kettenbruchentwicklung}, Nachr. Ges. Wiss. G\"ottingen Math-Phys.
Kl., 3, (1895), pp.~357-359.
\bibitem{Kle2}
F.~Klein, {\it Sur une repr\'esentation g\'eom\'etrique de
d\'eveloppement en fraction continue ordinaire}, Nouv. Ann. Math.
15(3), (1896), pp.~327--331.
\bibitem{Kon}
M.~L.~Kontsevich and Yu.~M.~Suhov, {\it Statistics of Klein
Polyhedra and Multidimensional Continued Fractions}, Amer. Math.
Soc. Transl., v.~197(2), (1999) pp.~9--27.
\bibitem{Kor0}
E.~I.~Korkina, {\it The simplest 2-dimensional continued
fraction}, International Geometrical Colloquium, Moscow 1993.
\bibitem{Kor1}
E.~I.~Korkina, {\it La p\'eriodicit\'e des fractions continues
multidimensionelles}, C. R. Ac. Sci. Paris, v. 319(1994),
pp.~777--780.
\bibitem{Kor2}
E.~I.~Korkina, {\it Two-dimensional continued fractions. The
simplest examples}, Proceedings of V.~A.~Steklov Math. Ins., v.
209(1995), pp. 143--166.
\bibitem{Kor3}
E.~I.~Korkina, {\it The simplest 2-dimensional continued
fraction.}, J. Math. Sci., 82(5), (1996), pp.~3680--3685.
\bibitem{Lac}
G.~Lachaud, {\it Poly\`edre d'Arnold et voile d'un c\^one
simplicial: analogues du th\`eoreme de Lagrange}, C. R. Ac. Sci.
Paris, v. 317(1993), pp.~711--716.
\bibitem{Lac2}
G.~Lachaud, {\it Voiles et Poly\`edres de Klein}, preprint n
95-22, Laboratoire de Math\'ematiques Discr\`etes du C.N.R.S.,
Luminy (1995).
\bibitem{Lac3}
G.~Lachaud, {\it Sails and Klein Polyhedra}, Contemp. Math.,
v.~210(1998), pp.~373-385.
\bibitem{LLL}
A.~K.~Lenstra, H.~W.~Lenstra, Jr., and L.~Lov\'asz, {\it
Factoring Polynomials with Rational Coefficients}, Mathematische
Ann., v.216(1982) pp.~515--534.
\bibitem{Min}
H.~Minkowski, {\it G\'en\'eralisation de le th\'eorie des
fractions continues}, Ann. Sci. Ec. Norm. Super. ser III, vol.
13, (1896), pp. 41--60.
\bibitem{Mou1}
J.-O.~Moussafir, {\it Sales and Hilbert bases}, Func. an. and
appl., v.~34(2000), n.~2, pp.~43--49.
\bibitem{Mou2}
J.-O.~Moussafir, {Voiles et Poly\'edres de Klein: Geometrie,
Algorithmes et Statistiques},
docteur en sciences th\'ese, Universit\'e Paris IX - Dauphine, (2000)\\
see also at http://www.ceremade.dauphine.fr/\~{}msfr/
\bibitem{Oka}
R.~Okazaki, {\it On an effective determination of a Shintani's
decomposition of the cone $\r^n_+$}, J. Math. Kyoto Univ.,
v33-4(1993), pp.~1057--1070.
\bibitem{Par1}
V.I.~Parusnikov, {\it Klein's polyhedra for the third extremal
ternary cubic form}, preprint 137 of Keldysh Institute of the
RAS, Moscow, (1995).
\bibitem{Par1.1}
V.I.~Parusnikov, {\it Klein's polyhedra for the fifth extremal
cubic form}, preprint 69 of Keldysh Institute of the RAS, Moscow,
(1998).
\bibitem{Par1.2}
V.I.~Parusnikov, {\it Klein's polyhedra for the seventh extremal
cubic form}, preprint 79 of Keldysh Institute of the RAS, Moscow,
(1999).
\bibitem{Par2}
V.I.~Parusnikov, {\it Klein's polyhedra for the fourth extremal
cubic form}, Mat. Zametki, 67(1), (2000), 110-128.
\bibitem{Shi}
T.~Shintani, {\it On evaluation of zeta functions of totally real
algebraic number fields at nonpositive integers}, J. Fac. Sci.
Univ. Tokyo Sect. IA, vol.~23(1976), pp.~393--417.
\bibitem{Sku1}
B.~F.~Skubenko, {\it Minima of a decomposible cubic form of three
variables}, Sci. Seminar Notes LOMI, vol.~168 (1988), Analytic
Number Theory and Theory of functions, 9, Leningrad, ''Nauka''.
\bibitem{Sku2}
B.~F.~Skubenko, {\it Minima of decomposible forms of degree $n$
of $n$ variables for $n \ge 3$}, Sci. Seminar Notes LOMI,
vol.~183 (1990), Modular functions and quadratic forms, 1,
Leningrad, ''Nauka''.
\bibitem{Tho}
E.~Thomas and A.~T.~Vasques, {\it On the resolution of cusp
singularities and the Shintani decomposition in totally real
cubic number fields}, Math. Ann. v.247(1980), pp.~1--20.
\bibitem{Tsu}
H.~Tsuchihashi, {\it Higher dimensional analogues of periodic
continued fractions and cusp singularities}, Tohoku Math. Journ.
v. 35(1983), pp.~176--193.
\bibitem{Voro}
G.~F.~Voronoi, {\it On one generalization of continued fraction
algorithm}, USSR Ac. Sci., v.1(1952),  pp.197-391.

\end{thebibliography}
\end{document}